\documentclass[JEDP]{cedram-sem}
\usepackage{amsfonts}
\usepackage{amsthm}

\newtheorem{theorem}{Theorem}[section]

\newtheorem{conjecture}[theorem]{Conjecture}

\usepackage{graphicx}
\newcommand{\R}{\mathbb{R}}
\newcommand{\N}{\mathbb{N}}
\DeclareMathOperator{\dive}{div}
\DeclareMathOperator{\curl}{curl}

\title 
[Navier-Stokes Controllability]
{On the controllability of the Navier-Stokes equation in a rectangle, with a little help of a distributed phantom force}

\author
[J.-M. \lastname{Coron}]
{\firstname{Jean-Michel} \lastname{Coron}}
\address{Laboratoire Jacques-Louis Lions, \\ Sorbonne Universit\'e, France}
\email{coron@ann.jussieu.fr}
\author
[F. \lastname{Marbach}]
{\firstname{Fr\'ed\'eric} \lastname{Marbach}}
\address{Univ Rennes, CNRS, France}
\email{frederic.marbach@ens-rennes.fr}
\author
[F. \lastname{Sueur}]
{\firstname{Franck} \lastname{Sueur}}
\address{Institut de Math\'ematiques de Bordeaux, \\ Universit\'e de Bordeaux, France}
\email{franck.sueur@math.u-bordeaux.fr}
\author
[P. \lastname{Zhang}]
{\firstname{Ping} \lastname{Zhang}}
\address{Academy of Mathematics \& Systems Science and Hua Loo-Keng Key Laboratory of Mathematic, The Chinese Academy of Sciences, China}
\email{zp@amss.ac.cn}


\begin{document}

\begin{abstract}
 This note echoes the talk given by the second author during the Journ\'ees EDP 2018 in Obernai. Its aim is to provide an overview and a sketch of proof of the result obtained by the authors in \cite{2018arXiv}, concerning the controllability of the Navier-Stokes equation. We refer the interested readers to the original paper for the full technical details of the proof, which will be omitted here, to focus on the main underlying ideas.
\end{abstract}

\maketitle

\section{Geometric setting}

\begin{figure}[ht!]
 \includegraphics{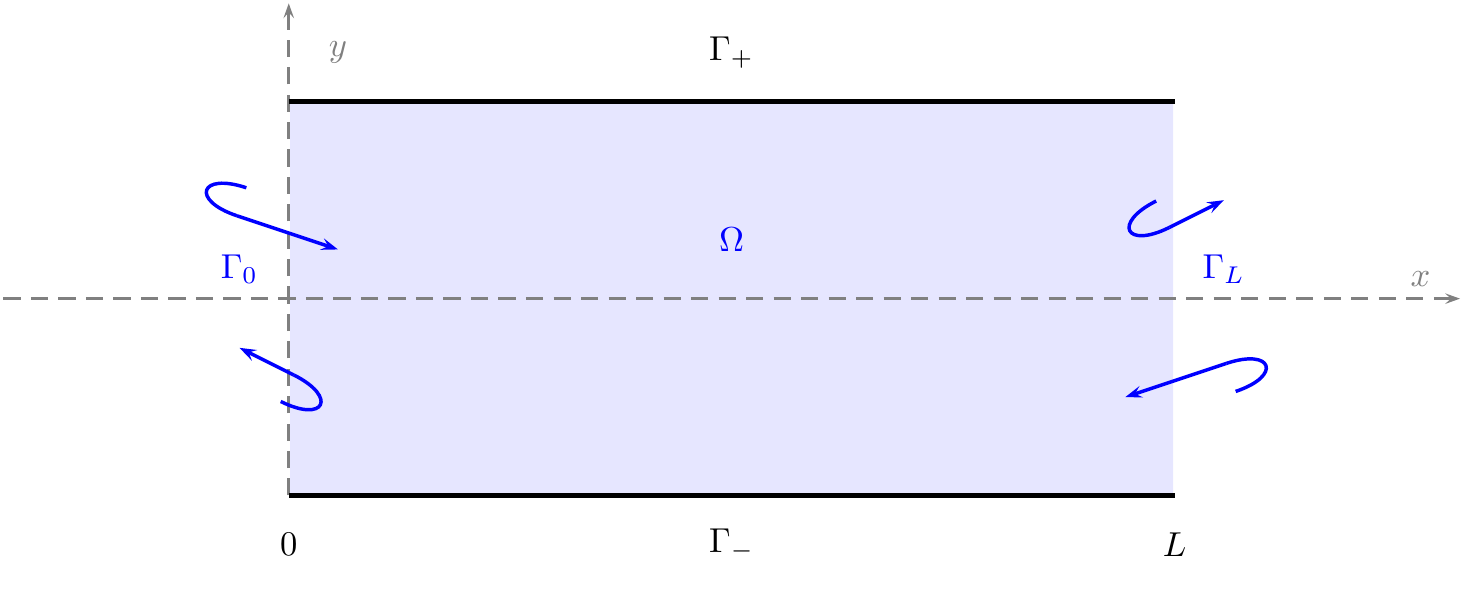}
 \caption{Physical domain $\Omega$}
 \label{fig:domain}
\end{figure}

We consider a rectangular domain $\Omega := (0,L) \times (-1, 1)$, where $L > 0$ is the horizontal length of the domain (see Figure \ref{fig:domain}). We will use $(x,y) \in \Omega$ as coordinates. We see this rectangular domain as a tube or a river, in the interior of which a fluid evolves. During some time interval $[0,T]$, the evolution of the fluid velocity $u(t,x,y)$ is governed by the homogeneous incompressible Navier-Stokes equation:
\begin{equation} \label{ns}
 \left\{
  \begin{aligned}
   \partial_t u + (u \cdot \nabla) u - \Delta u + \nabla p = f, \\
   \dive u = 0,
  \end{aligned}
 \right.
\end{equation}
where $f(t,x,y)$ is a small external vectorial forcing term, whose role will be explained below and $p(t,x,y)$ is the scalar pressure field corresponding to the incompressibility constraint.

On the upper and lower horizontal boundaries $\Gamma_\pm := (0,L) \times \{ \pm 1 \}$, corresponding to the walls of the tube or the banks of the river, we assume that the fluid satisfies the usual no-slip Dirichlet boundary condition:
\begin{equation} \label{dirichlet}
 u = 0.
\end{equation}
Conversely, a key feature of the geometric setting at stake is that no boundary condition is prescribed \emph{a priori} on the left and right vertical boundaries $\Gamma_0 := \{0\} \times (-1,1)$ and $\Gamma_L := \{L\} \times (-1,1)$. This under-determination models the idea that we can act on the system by exerting some forcing (say, through suction or blowing actions) on the fluid.

\section{Cauchy problem}

In 2D, it is known that weak Leray solutions to the homogeneous incompressible Navier-Stokes equation exist globally and are unique. In our setting, uniqueness is not guaranteed because the problem is under-determined due to the possible choices of the boundary conditions on $\Gamma_0$ and $\Gamma_L$ (which correspond to controls).

More precisely, let $L^2_{\dive}(\Omega)$ denote the space of $L^2$ vector fields on $\Omega$ which are divergence-free and tangent to the boundaries $\Gamma_\pm$. Given $T > 0$, an initial data $u_* \in L^2_{\dive}(\Omega)$, and a forcing $f \in L^1((0,T);L^2(\Omega))$, we will say that $u \in C^0([0,T];L^2_{\dive}(\Omega)) \cap L^2((0,T);H^1(\Omega))$ is a weak Leray solution to \eqref{ns} and \eqref{dirichlet} with final data $u_T \in L^2_{\dive}(\Omega)$ when it satisfies the weak formulation:
\begin{equation}
 \begin{split}
  - \int_0^T \int_\Omega u \cdot \partial_t \varphi 
  + \int_0^T \int_\Omega & (u\cdot\nabla) u \cdot \varphi
  + 2 \int_0^T \int_\Omega D(u) : D(\varphi) \\
  & = \int_\Omega u_* \cdot \varphi(0,\cdot) 
  - \int_\Omega u_T \cdot \varphi(T,\cdot)
  + \int_0^T \int_\Omega \varphi \cdot f,
 \end{split}
\end{equation}
for every test function $\varphi \in C^\infty([0,T]\times\bar{\Omega})$ which is divergence-free, tangent to $\Gamma_\pm$ and vanishes on $\Gamma_0$ and $\Gamma_L$.

Another way to formulate the Cauchy problem is to see weak Leray solutions on $\Omega$ as the restriction to the physical domain $\Omega$ of weak solutions defined on a larger domain, say the strip $\mathfrak{B} := \R \times (-1,1)$, corresponding to some extensions of the initial data and of the external force. Given any (reasonable) choice of extensions for $u_*$ and $f$, there exists a unique global weak solution on $\mathfrak{B}$, which can then be restricted to $\Omega$.

\section{A conjecture of Lions}

In the late 1980's, Jacques-Louis Lions formulated multiple open problems and conjectures concerning the controllability of systems governed by partial differential equations. In particular, in \cite{MR1147191}, he asked whether the Navier-Stokes equation was small-time globally null controllable. There are many ways to set this question, depending on the geometry, on the exact goals, and on the nature of the exerted controls (which can either be a distributed force in some strict subset of the domain or come into play through boundary data). In our geometrical setting, the conjecture of controllability can be formulated as:

\begin{conjecture}
 Let $T > 0$ and $u_* \in L^2_{\dive}(\Omega)$. There exists a weak Leray solution to~\eqref{ns} with $f = 0$ and \eqref{dirichlet} such that the final state satisfies $u(T,\cdot) = 0$.
\end{conjecture}

The difficulty in the question comes from the combination of multiple factors. First, the allotted control time $T > 0$ may be very small, which requires to use an asymptotically rapid strategy. Second, the initial data $u_*$ may be very large, so that the nonlinearity in the Navier-Stokes equation plays an important role. Last, but not least, the controls are only exerted on a strict subset $\Gamma_0 \cup \Gamma_L$ of the full boundary $\partial\Omega$. One can expect that specific phenomenons occur near the uncontrolled parts $\Gamma_\pm$.

\section{Our main controllability result}

In \cite{2018arXiv}, we proved a result which almost brings a positive answer to the above conjecture. Whereas the initial conjecture implies to find an exact solution of the Navier-Stokes equation with a null forcing term, we introduce a non-zero but arbitrarily small forcing, in arbitrarily strong norms.

\begin{theorem} \label{thm}
 Let $T > 0$ and $u_* \in L^2_{\dive}(\Omega)$. For every $k \in \N$ and every $\eta > 0$, there exists a force $f \in L^1((0,T); H^k(\Omega))$ satisfying
 \begin{equation} \label{phantomatique}
  \| f \|_{L^1((0,T);H^k(\Omega))} \leq \eta,
 \end{equation}
 and an associated weak Leray solution $u \in C^0([0,T];L^2_{\dive}(\Omega)) \cap L^2((0,T);H^1(\Omega))$ to \eqref{ns} and \eqref{dirichlet} satisfying $u(0) = u_*$ and $u(T) = 0$.
\end{theorem}

In this under-determined formulation of the control result, the boundary controls (i.e.\  the traces of $u$ on the boundaries $\Gamma_0$ and $\Gamma_L$) are not explicitly written.

Hence, we almost obtain small-time global exact null controllability. Our method does
not easily extend to obtain the ``true'' control result with $f = 0$. Indeed, one cannot pass to the limit in the main theorem because there is no \emph{a priori} bound on the size of the trajectories $u$ as $\eta \to 0$.

The small correction we need is linked with our proof strategy (which creates a boundary layer) and our proof technique (which relies on horizontal analyticity). It is likely that proving the result for $f = 0$ requires both a new strategy and a new technique.

The fact that $\Omega = (0,L) \times (-1,1)$ is a ``flat'' domain is also very important for our proof. More precisely, the key point is that the uncontrolled boundaries $\Gamma_\pm$ are flat in the horizontal direction. This feature allows us to introduce almost explicit expressions for some of the profiles that build up the solution $u$, which are solutions to linear equations.

\section{Discussion on earlier results} \label{sec:earlier}

The open problems introduced by Jacques-Louis Lions concerning controllability for fluid mechanics problems have received a large attention.

\paragraph{Small initial data and local results}

Small-time local null controllability was already known. For every $T > 0$, there exists $\delta_T > 0$ such that, for every $u_* \in L^2(\Omega)$ satisfying $\|u_*\| \leq \delta_T$, one can find controls driving $u_*$ to the null equilibrium state $u = 0$ in time $T$. This can be done using only boundary controls, without any distributed force ($f = 0$). In this case, since the state is small, one sees the bilinear term in the Navier-Stokes system as a small perturbation term of the Stokes equation so that the controllability is proved thanks to Carleman estimates and fixed point theorems. Loosely speaking, such an approach corresponds to low Reynolds controllability. We refer to \cite{MR2103189,MR1308746,MR1804497} for some important contributions to this topic, successively improving the smallness assumptions, the control domains or the reachable targets.
 
\paragraph{Global results without boundaries} 
 
For large initial data, a setting corresponding to controllability at large Reynolds numbers, the first author and Fursikov proved global null controllability for the Navier-Stokes system in a 2D manifold without boundary in~\cite{MR1470445} (in this case, the control is an internal control exerted from a small open subset of the domain). In~\cite{MR1728643}, Fursikov and Imanuvilov proved a small-time global control result when the control is exerted on the full boundary $\partial\Omega$ of the physical fluid domain. Both geometries share the important feature that there is no uncontrolled portion of the boundary.

\paragraph{Navier slip-with-friction boundary condition} Jacques-Louis Lions' problem has been solved in~\cite{2016arXiv161208087C} by the first three authors in the particular case of the Navier slip-with-friction boundary condition (see also \cite{2017arXiv170307265C} for a gentle introduction to this result). This boundary condition is less stringent than \eqref{dirichlet} since it allows the fluid to slide tangentially along the boundary. In this context, small-time global exact null controllability and small-time global exact controllability to trajectories hold for every regular domain (2D and 3D) and for every subset of the domain where the control is exerted, provided that it intersects each connected component of the boundary of the physical domain.

\paragraph{Partial results with large forcing}

The closest works to Theorem \ref{thm} are references~\cite{MR2269867, MR2994698}, in which related results are obtained in very similar settings. 
These works prove a version of Theorem \ref{thm} for which the distributed force $f$ can be chosen small in $L^p((0,T);H^{-1}(\Omega))$, where $1 < p < 4/3$.
The fact that our phantom force can be chosen arbitrarily small in the space $L^1((0,T),H^k(\Omega))$ for any $k \geq 0$, is the major improvement of this work. In particular, being small, say in $C^1(\bar{\Omega})$ guarantees that there is not fast scale variations of our distributed force near the uncontrolled boundaries. This possibility is not ruled out by a conclusion on the smallness in $H^{-1}(\Omega)$ of the forcing term.

\section{A strategy based on the flushing of the vortexes}

If one thinks that the vector field $u(t,\cdot)$ is described by the combination of its potential part and its vorticity, driving to zero requires to drive both parts to zero. Thanks to the incompressibility constraint, it is very easy to make the potential part vanish, almost instantly. Indeed, if one chooses null boundary controls on $\Gamma_0$ and $\Gamma_L$, then at any instant $t > 0$, the full state $u(t,\cdot)$ can be recovered from its vorticity $\omega(t,\cdot)$ through the following div-curl problem:
\begin{equation}
 \left\{
  \begin{aligned}
   \curl u(t,\cdot) & = \omega(t,\cdot) 
   && \text{in } \Omega, \\
   \dive u(t,\cdot) & = 0
   && \text{in } \Omega, \\
   u(t,\cdot) \cdot n & = 0
   && \text{on } \partial\Omega.
  \end{aligned}
 \right.
\end{equation}
We can thus assume that the initial data has a vanishing average horizontal velocity, i.e.\ $\int_\Omega u_* \cdot e_x = 0$, where $e_x$ is the tangential unit vector. If it is not the case, using such null controls will ensure it for any positive time. 

We embed $\Omega$ in the band $\mathfrak{B} = \R \times (-1,+1)$ and extend the initial data $u_*$ to a compactly supported (say on $[-L,2L] \times [-1,1]$) divergence-free initial data on $\mathfrak{B}$ (this is possible when $u_*$ has zero average tangential speed), which we will still denote by $u_*$. We work in the extended domain $\mathfrak{B}$ for simplicity.

Our goal is thus to build a solution such that $u(T)_{\rvert \Omega} = 0$. In fact, it is sufficient to achieve $\| u(T)_{\rvert \Omega} \|_{L^2(\Omega)} \ll 1$, since local controllability is known for the Navier-Stokes equation (see the paragraph \emph{Small initial data and local results} of Section \ref{sec:earlier}).

\begin{figure}[ht!]
 \includegraphics{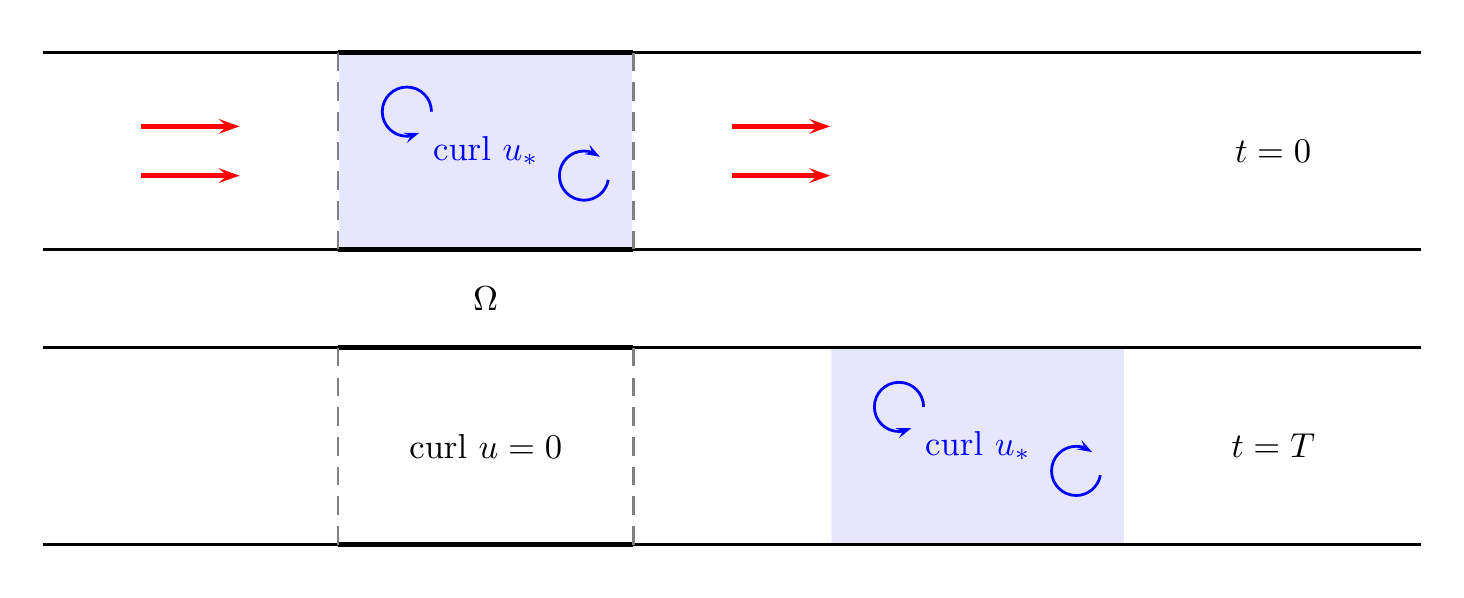}
 \caption{Flushing process for the vorticity}
 \label{fig:flush}
\end{figure}

Recalling that, in 2D, the vorticity is transported by the flow, the first important idea is to flush the support of the initial vorticity $\omega_* := \curl u_*$ outside of the initial physical domain $\Omega$, into the extension $\mathfrak{B}\setminus\Omega$. We perform this task using the incompressibility and introducing artificially a high pressure gradient as sketched in Figure \ref{fig:flush}. 

\section{Asymptotic implementation of the flushing method}

In order to implement the intuition of Figure \ref{fig:flush}, we introduce a small parameter $\varepsilon > 0$ and we will construct a solution $u(t,x,y)$ given under the form
\begin{equation}
 u(t,x,y) = \frac{1}{\varepsilon} u^{\varepsilon}\left(\frac{t}{\varepsilon}, x, y\right),
\end{equation}
where the new unknown $u^\varepsilon$ must now solve the following modified equation on a larger time interval $t \in (0,T/\varepsilon)$
\begin{equation} \label{ns.eps}
 \left\{
\begin{aligned}
\partial_t u^\varepsilon + (u^\varepsilon \cdot \nabla) u^\varepsilon - \varepsilon \Delta u^\varepsilon + \nabla p^\varepsilon & = f^\varepsilon, \\
\dive u^\varepsilon & = 0, \\
u^\varepsilon(0) & = \varepsilon u_*,
\end{aligned}
\right.
\end{equation}
where we introduced $p^\varepsilon(t,x,y) = \varepsilon^2 p(\varepsilon t, x, y)$ and $f^\varepsilon(t,x,y) = \varepsilon^2 f(\varepsilon t, x, y)$. Within this scaling, the goal is to construct a solution such that $\| u^\varepsilon(T/\varepsilon)_{\vert \Omega} \|_{L^2(\Omega)} \ll \varepsilon$.

Heuristically, we wish to build a solution to \eqref{ns.eps} which behaves as
\begin{equation} \label{exp.1}
 u^\varepsilon(t,x,y) \approx h(t) e_x + \varepsilon u^1(t,x,y) + o(\varepsilon),
\end{equation}
where $h \in C^\infty(\R_+;\R)$ is supported on $(0,T)$ and has a sufficiently large integral, say $\int_0^T h(t) \mathrm{d} t \geq 3 L$, and $u^1$ is the solution to the linearized version of \eqref{ns.eps} around $h(t)e_x$ (which is a solution of the underlying Euler equation, see the red arrows on Figure \ref{fig:flush}),
\begin{equation} \label{u1}
\left\{
\begin{aligned}
\partial_t u^1 + h \partial_x u^1 + \nabla p^1 & = 0, \\
\dive u^1 & = 0, \\
u^1(0) & = u_*.
\end{aligned}
\right.
\end{equation}
Of course, thanks to the simple geometrical setting, \eqref{u1} can be solved explicitly as
\begin{equation}
 u^1(t,x,y) = u_* \left( x - \int_0^t h(t') \mathrm{d} t', y \right).
\end{equation}
In particular, if $u_*$ was compactly supported, say on $[-L,2L] \times [-1,1] \subset \bar{\mathfrak{B}}$, then $u^1$ vanishes inside $\Omega$ for $t \geq T$ thanks to the assumption that $\int h \geq 3L$.

If we believe that the remainder in \eqref{exp.1} is indeed $o(\varepsilon)$, then the theorem is proved since, for $t \geq T$ (thus including $t = T/\varepsilon$), $h(t)$ vanishes and $u^1(t)$ vanishes inside $\Omega$, so that $\| u^\varepsilon(T/\varepsilon)_{\vert \Omega} \|_{L^2(\Omega)} \ll \varepsilon$.

\section{Tangential boundary layers}

Unfortunately, the leading order profile $h(t) e_x$ is the solution of the underlying Euler equation (corresponding to $\varepsilon = 0$ in \eqref{ns.eps}) and only satisfies the normal impermeability condition $u \cdot e_y = 0$ on $\Gamma_\pm$. The tangential boundary condition $u \cdot e_x = 0$ is not satisfied by this profile. Hence, there is no chance for an expansion like \eqref{exp.1} to hold.

This discrepancy is very usual when studying the convergence of Navier-Stokes to Euler in the vanishing viscosity setting. It gives rise to the theory of boundary layers: a small region, here of width $\varepsilon^{\frac12}$, within which the viscous effects remain important and allowing to recover the missing boundary condition. Plugging such an \emph{Ansatz} depending on a fast variable in the Navier-Stokes equations yields the Prandtl equation \cite{prandtl1904uber} governing the evolution of the boundary layer profile. Here, thanks to the flat geometric setting and the invariance with respect to $x$ of the main profile, they take a particularly simple form. Indeed, we change our expansion \eqref{exp.1} into
\begin{equation} \label{exp.2}
 u^\varepsilon(t,x,y) \approx \left[ h(t) - V\left(t,\frac{1 + y}{\sqrt{\varepsilon}}\right) \right] e_x
 + \varepsilon u^1(t,x,y) + o(\varepsilon),
\end{equation}
where $V : \R_+ \times \R_+ \to \R$ is the solution to the following heat equation (a very simplified version of the Prandtl equation in our setting):
\begin{equation} \label{V}
\left\{
\begin{aligned}
\partial_t V + \partial_{zz} V & = 0, \\
V(t,0) & = h(t), \\
V(0,z) & = 0.
\end{aligned}
\right.
\end{equation}
In fact, a second symmetrical corrector depending on $1-y$ is required in order to account for the boundary layer near $\Gamma_+$, and smooth slowly varying cutoff functions are needed in order to avoid interaction between the two correctors. We will not consider these details here and proceed with the computations only with the corrector near $\Gamma_-$, as they already contain the core ideas.

These correctors allow to build a reference flow which fully satisfies the boundary conditions on $\Gamma_\pm$, enabling us to hope to prove \eqref{exp.2}.

\section{Main difficulties} \label{sec:difficulties}

Although the boundary correctors only change the value of the reference flow in small strips near the boundaries, they introduce two important difficulties with respect to our controllability goal.

First, at the final time $t = T/\varepsilon$, although $h(t)$ vanishes and $u^1(t)$ vanishes inside $\Omega$, it is not the case for $V$. More precisely, one has
\begin{equation} \label{vlarge}
 \left\| u(T)_{\rvert \Omega} \right\|_{L^2(\Omega)}
 = \left\| \frac{1}{\varepsilon} u^\varepsilon\left(\frac{T}{\varepsilon}\right)_{\rvert \Omega} \right\|_{L^2(\Omega)}
 \approx \varepsilon^{-\frac34}
 \left\| V\left(\frac{T}{\varepsilon}\right) \right\|_{L^2(\R_+)}.
\end{equation}
For $t \geq T$, the heat equation \eqref{V} has zero source term and the profile $V$ decays in $L^2(\R_+)$. Unfortunately, without any additional assumption, studying the decay rates for the free heat equation on the half line only yields a weak decay of the form $\| V(t) \|_{L^2(\R_+)} \leq C t^{-\frac 14}$ as $t \to +\infty$, which is not sufficient to counterbalance the prefactor of \eqref{vlarge}.

Second, trying to make expansion \eqref{exp.2} rigorous and computing the equation satisfied by the remainder $r^\varepsilon$ (for $u^\varepsilon = h + V + \varepsilon u^1 + \varepsilon r^\varepsilon$), yields an evolution equation with a bad amplification term:
\begin{equation} \label{r}
\left\{
\begin{aligned}
\partial_t r^\varepsilon 
+ \varepsilon^{-\frac12} r^\varepsilon_2 \partial_z V e_x
+ A_\varepsilon r^\varepsilon
+ \varepsilon (r^\varepsilon \cdot \nabla) r^\varepsilon
- \varepsilon \Delta r^\varepsilon
 & = \sigma^\varepsilon, \\
\dive r^\varepsilon & = 0, \\
r^\varepsilon(0) & = 0.
\end{aligned}
\right.
\end{equation}
In \eqref{r}, $\sigma^\varepsilon$ is a small source term in some appropriate sense (one can think $\sigma^\varepsilon = o(1)$ in $L^1((0,T/\varepsilon);L^2(\mathfrak{B}))$ for example). The amplification has a reasonable part $A_\varepsilon r^\varepsilon$ (one can think that its norm in $L^1((0,T/\varepsilon); L^\infty(\mathfrak{B}))$ is bounded uniformly with respect to $\varepsilon$) and a very bad part $\varepsilon^{-\frac12} r^\varepsilon_2 \partial_z V e_x$.
Performing naive Gr\"onwall estimates on this equation is therefore bound to fail due to this term, even more so since we intend to perform these estimates up to the large final time $T/\varepsilon$ (we would then expect an exponential amplification of the form $\exp(\varepsilon^{-3/2})$). So, \emph{a priori}, the remainder is not small.

\section{Recasting amplification as a loss of derivative}

We start by dealing with the second problem. Using the divergence free condition on $r^\varepsilon$ and the null boundary condition on $r^\varepsilon_2$, we wish to rewrite the amplification term. We perform the computation near the lower wall $y=-1$. Near this wall, $V$ is evaluated at $z = \varepsilon^{-\frac 12} (1 + y)$. Hence one has
\begin{equation}
 \begin{split}
  \varepsilon^{-\frac12} r^\varepsilon_2(t,x,y) \partial_z V(t,z) e_x
  & = \varepsilon^{-\frac12}(1+y) \left( \frac{1}{1+y} \int_{-1}^y \partial_y r^\varepsilon_2(t,x,y') \mathrm{d} y' \right) \partial_z V(t,z) e_x \\
  & = - \left(\frac{1}{1+y} \int_{-1}^y \partial_x r^\varepsilon_1(t,x,y') \mathrm{d} y' \right)
   z \partial_z V(t,z) e_x.
 \end{split}
\end{equation}
Thus, the amplification term has been recast as a local average in the normal direction of $\partial_x r^\varepsilon_1$. However, this term does not have the structure of a transport term that would disappear during energy estimates by integration by parts. On a formal level, one should rather think of this term as the structure-less term
\begin{equation} \label{loss}
  (z \partial_z V) |\partial_x| r^\varepsilon,
\end{equation}
where $|\partial_x|$ is defined as the Fourier multiplier by $|\xi|$, where $\xi$ is the horizontal Fourier variable. Since $V$ is the solution to \eqref{V}, for each $t \geq 0$, the map $z \mapsto z \partial_z V(t,z)$ belongs to $L^\infty(\R_+)$ because $V(t,\cdot)$ and its derivatives decay exponentially (with respect to $z \to +\infty$).

There is \emph{a priori} no hope to ``absorb'' a term such as \eqref{loss} by the $- \varepsilon \Delta r^\varepsilon$ dissipation term of \eqref{r}, because the estimate would once again degenerate as $\varepsilon \to 0$. Instead, we think of \eqref{loss} as a loss of derivative, and we will work in an analytic setting (with respect to the tangential variable), so that loosing one derivative (among an infinite number of derivatives) is not too bad.

In the context of Navier-Stokes boundary layers, analyticity was first used in \cite{MR1617542,MR1617538} to prove both the existence of solutions to the Prandtl equation and the convergence of the vanishing viscosity Navier-Stokes solution to an Euler+Prandtl system for analytic data.

\section{Cauchy-Kowaleskaya schema}

Due to the term \eqref{loss}, the analytic radius of the solution $r^\varepsilon$ will decay as time increases. This rough idea can be very precisely quantified thanks to an idea linked with Cauchy-Kowaleskaya type theorems. Let $\rho \in C^1(\R_+;\R)$. We introduce the new unknown 
\begin{equation} \label{ck.trick}
 r^\varepsilon_\rho := e^{\rho(t) |\partial_x|} r^\varepsilon.
\end{equation}
This change of unknown is licit for example when the tangential Fourier transform of $r^\varepsilon$ is supported on some bounded region $-N \leq \xi \leq N$, so this trick has to be performed on a ``frequency-truncated'' version of \eqref{r}, which will then pass to the limit since the resulting estimates will not depend on $N$.
Under the change of unknown \eqref{ck.trick}, equation \eqref{r} is roughly changed into
\begin{equation} \label{ck.balance}
 \partial_t r^\varepsilon_\rho - \rho'(t) |\partial_x| r^\varepsilon_\rho 
 - (z \partial_z V) |\partial_x | r^\varepsilon_\rho = ...
\end{equation}
Therefore, multiplying \eqref{ck.balance} by $r^\varepsilon_\rho$ and using Parseval's formula yields
\begin{equation}
 \frac 12 \frac{\mathrm{d}}{\mathrm{d}t} \int (\hat{r}^\varepsilon_\rho)^2
 - \rho'(t) \int |\xi| (\hat{r}^\varepsilon_\rho)^2 
 \leq \| z \partial_z V(t,z) \|_{L^\infty(\R_+)} \int |\xi| (\hat{r}^\varepsilon_\rho)^2  + ...,
\end{equation}
so that the derivative loss term of the right-hand side can be absorbed if and only if
\begin{equation} \label{rho'}
 - \rho'(t) \geq \| z \partial_z V(t,z) \|_{L^\infty(\R_+)}.
\end{equation}
Since \eqref{rho'} must be satisfied for $t \in [0,T/\varepsilon]$ and since we wish $\rho$ to stay positive, we need to choose an initial analyticity radius $\rho(0)$ such that
\begin{equation} \label{int.large}
 \rho(0) \geq \int_0^{+\infty} \| z \partial_z V(t,z) \|_{L^\infty(\R_+)} \mathrm{d}t.
\end{equation}
\emph{A priori}, there is no reason for this integral to be finite, so we will need to adapt our construction in order to ensure it.

\section{Preparing a good-enough dissipation}

We now turn to the first problem mentioned in Section \ref{sec:difficulties}: namely the fact that the boundary layer term $V$ is not small enough at the final time $T/\varepsilon$. We wish to choose the source term $h$ of \eqref{V} more wisely in order to ensure that $V$ decays sufficiently fast. As an added benefit, this will make the integral in \eqref{int.large} finite.

For $t\geq T$, $h(t) = 0$ so \eqref{V} is a free heat equation with null boundary condition at $z = 0$. The decay rate of the free heat equation on the half-line is linked with the low frequencies of the ``initial'' data $V(T,\cdot)$. More precisely, it depends on the number of vanishing derivatives of its Fourier transform at zero. These quantities are linked to the $z$-moments $\int_0^{+\infty} z^k V(T,z) \mathrm{d}z$, which are linked with the $t$-moments of the source term $h$, $\int_0^T t^k h(t) \mathrm{d} t$. For example, choosing $h \in C^\infty([0,T];\R)$ such that $\int_0^T h(t) \mathrm{d} t = 0$ guarantees that $\int_0^{+\infty} z V(T,z) \mathrm{d}z  = 0$, which in turn improves the decay rate of the solution by a factor $1/t$ for $t \to +\infty$.

The key idea here is thus to choose a function $h \in C^\infty([0,T];\R)$ which has a finite number of null time moments. This guarantees that the solution to \eqref{V} will decay sufficiently fast (not only in $L^2(\R_+)$ but also for stronger functional spaces including polynomial weights in $z$ and Sobolev norms). As a consequence, we obtain
\begin{equation}
 \left\| V \left(\frac T{\varepsilon}\right) \right\|_{L^2(\R_+)} = O(\varepsilon^3) 
 \quad \textrm{and} \quad
 \int_0^{+\infty} \| z \partial_z V(t,z) \|_{L^\infty(\R_+)} \mathrm{d}t
 < + \infty. 
\end{equation}

\section{Killing the initial data when it is outside}

The initial intuition, depicted in Figure \ref{fig:flush} was to choose $\int_0^T h(t) \mathrm{d} t \geq 3L$ in order to flush the initial vorticity $\omega_*$ outside of the physical domain $\Omega$. Now that we need to choose $\int_0^T h(t) \mathrm{d}t = 0$, this intuition is not sufficient anymore. However, we can use controls to suppress the initial vorticity while it is outside of the physical domain.

For example, is we choose $h$ such that $\int_0^{T/3} h(t) \mathrm{d}t = 3L$, $h = 0$ on $(T/3,2T/3)$ and $\int_{2T/3}^{T} h(t) \mathrm{d}t = -3L$, we have a reference flow which is globally of zero average, but for which there exists an intermediate time when the initial vorticity $\omega_*$ is fully outside of the physical domain.

During $(T/3,2T/3)$, we thus apply a control (in the form of a source term in \eqref{u1}, supported outside of $\Omega$), which is designed to obtain $u^1(2T/3) = 0$. Hence, when $h$ becomes negative and ``brings back'' fluid particles into $\Omega$, it carries only a vanishing vorticity.

Heuristically, one sets
\begin{equation} \label{u1.beta}
 u^1(t,x,y) = \beta(t) u_* \left( x - \int_0^t h(t') \mathrm{d}t', y \right),
\end{equation}
where $\beta \in C^\infty(\R_+;[0,1])$ is such that $\beta(t) = 1$ for $t \leq T/3$ and $\beta(t) = 0$ for $t \geq 2T/3$. This defines a solution of \eqref{u1} with a non-zero right-hand side $f^1$ supported in $[T/3,2T/3] \times [2L,5L] \times [-1,1]$:
\begin{equation} \label{f1}
 f^1(t,x,y) := \beta'(t) u_* \left( x - \int_0^t h(t') \mathrm{d}t', y \right).
\end{equation}
Technically, one should write a formula like \eqref{u1.beta} on the stream function in order to preserve the divergence-free condition.

\section{Dealing with the non-linearity with Chemin's method}

An important drawback of the change of unknown \eqref{ck.trick} is that it destroys the nice structure of the nonlinear term $\varepsilon (r^\varepsilon \cdot \nabla) r^\varepsilon$ in equation \eqref{r}. Usually, this term disappears during the standard $L^2$ energy estimate obtained by multiplying equation \eqref{r} by $r^\varepsilon$ and using the divergence-free condition. This simplification does not happen anymore after our change of unknown and we must estimate this term.

Using an idea introduced by Chemin in \cite{MR2145938}, we now see $\rho$ as the unknown solution to the highly nonlinear ODE
\begin{equation} \label{ode}
 \rho'(t) = - \| z \partial_z V(t,z) \|_{L^\infty(\R_+)} - \varepsilon \| \nabla r^\varepsilon_\rho(t) \|_{\dot{B}^0_{2,1}},
\end{equation}
where $\dot{B}^0_{2,1}$ is a homogeneous Besov space associated with frequency truncations in the tangential direction and is designed to have the critical Sobolev embedding in two dimensions $(r^\varepsilon_\rho, \nabla r^\varepsilon_\rho) \in \dot{B}^0_{2,1} \Rightarrow r^\varepsilon_\rho \in L^\infty$. Exploiting the divergence-free condition on $r^\varepsilon_\rho$ and this new definition of $\rho'$ allows to control the nonlinear term.

However, since \eqref{ode} is a nonlinear ODE, it is not clear \emph{a priori} that its solution stays well defined (and positive) up to the final time $T/\varepsilon$ (one could have $\rho \to -\infty$ in finite time and we need to ensure $\rho > 0$ to stay within the analytic setting). Hence, we must perform a parallel estimate for $\rho$ in the same time as we are estimating $r^\varepsilon_\rho$. Here, we use the viscous term $-\varepsilon \Delta r^\varepsilon$ of equation \eqref{r}. Indeed, by Cauchy-Schwarz, the total decay of $\rho$ can be bounded as
\begin{equation}
 \varepsilon \int_0^{T/\varepsilon} \| \nabla r^\varepsilon_\rho \|_{\dot{B}^0_{2,1}}
 \leq \sqrt{T} \left( \varepsilon \int_0^{T/\varepsilon} \| \nabla r^\varepsilon_\rho \|_{\dot{B}^0_{2,1}}^2 \right)^{\frac 12},
\end{equation}
and the right-hand side is precisely the type of quantity for which we obtain bounds thanks to the viscous term $-\varepsilon \Delta r^\varepsilon$ when we perform $\dot{B}^0_{2,1}$ energy estimates on \eqref{r}.

Working with $\dot{B}^0_{2,1}$ (rather than $L^2$) is necessary in our context to benefit from the embedding mentioned above. The nonlinear term $r \nabla r$ is then estimated using paradifferential calculus techniques (including Bony's paraproducts), notably inspired by
\cite{MR2776367, MR3464051}.

\section{Spotting the phantoms}

A drawback of the analytic setting considered above is the use of a phantom force (in the sense of a source term supported everywhere, arbitrarily small in an arbitrarily strong Sobolev space) for two different purposes, which we reveal here. Of course, our strategy also requires large source terms (the controls) which are exclusively supported outside of the physical domain $\Omega$.

\paragraph{Analytic regularization of the initial data} First, we need the initial data $u_*$ to be analytic. Since Theorem \ref{thm} is stated with an $L^2$ initial data, we need a strategy to regularize it. It is well known that the Navier-Stokes equation exhibits a strong smoothing effect (thanks to the dissipation term) and that the solution instantly becomes analytic. However, the analytic radius at time $t > 0$ is only known to grow like $\sqrt{t}$. Since we seek a small-time control result, the natural smoothing only yields a small analyticity radius. However, the total loss of analytic radius in our setting is linked with the quantity
\begin{equation} \label{total.loss}
 \int_0^{+\infty} \| z \partial_z V(t,z) \|_{L^\infty(\R_+)} \mathrm{d} t.
\end{equation}
In turn, this quantity depends on $L$ and $T$ through the choice of the base flow $h$. It can be checked that since we require $\int_0^{T/3} h \geq 3L$, the quantity \eqref{total.loss} is bounded below. Hence, as a first step of our result, we use an external source term supported everywhere to trim off the high tangential frequencies of the initial data and make it analytic with a sufficient radius (say twice the value of \eqref{total.loss}). 

Since our method only needs to know that the analytic radius is large enough (and not that the associated analytic norm of the initial data is small), this clipping process can be done with a small source term even in a strong Sobolev space, ensuring~\eqref{phantomatique}.

\paragraph{Almost compactly supported extension} Second, looking at \eqref{f1} defining the external force used to drive $u^1$ to zero, one sees that its size within $\Omega$ is linked to the values of the extension $u^*$ in $[-3L,-2L] \times [-1,1]$. Our initial idea was to choose $u^*$ compactly supported, say in $[-L,2L] \times [-1,1]$. Of course, since we need $u_*$ to be analytic in the tangential direction, it cannot simultaneously have a compact support in $x$. The most we can require is that the Sobolev norm of the analytic extension $u_*$ is small in $[-3L,-2L] \times [-1,1]$. Then, from \eqref{f1}, we see that $f^1$ can be split as a control part (large, but supported outside of $\Omega$) and a phantom part (small, but supported inside $\Omega$). 

In fact, our detailed construction also proves that we can localize the support of this second phantom force in the vertical direction so that it does not touch the horizontal boundaries~$\Gamma_\pm$. More precisely, for every $T > 0$ and $u^* \in L^2(\Omega)$, we prove that there exists $\delta > 0$ such that, for any $k \in \N$ and $\eta > 0$, we can maintain the result of Theorem~\ref{thm} while ensuring that $\mathrm{supp} f^1_{\rvert \Omega} \subset [0,L] \times [-1+\delta,1-\delta]$. This highlights the fact that the main role of this second phantom force is to allow us to work in an analytic setting (but not to take care directly of the boundary layer).

\def\cprime{$'$}

\end{document}